\documentclass{elsart}

 \usepackage{epsfig}

\usepackage{amssymb}

\def\D{\mathcal{D}}
\def\H{\mathcal{H}}

\begin{document}

\begin{frontmatter}



\title{Improved upper bounds for vertex and edge fault diameters of Cartesian graph bundles}



\author[label2,label3]{Rija Erve\v s}, 
\ead{rija.erves@um.si}
\author[label2,label4]{Janez \v{Z}erovnik}
\ead{janez.zerovnik@imfm.uni-lj.si}

\address[label2]{
Institute of Mathematics, Physics and Mechanics,\\ Jadranska 19, Ljubljana 1000, Slovenia.\\}
\address[label3]{
FCE,
University of Maribor,\\ Smetanova 17, Maribor 2000, Slovenia.\\}
\address[label4]{
FME,
University of Ljubljana,\\ A\v{s}ker\v{c}eva 6, Ljubljana 1000, Slovenia.\\}

\begin{abstract} 
Mixed fault diameter of a graph $G$, $ \D_{(a,b)}(G)$, is the maximal diameter of $G$ after deletion of any $a$ vertices and any $b$ edges.
Special cases are the (vertex) fault diameter $\D^V_{a} =  \D_{(a,0)}$ and the edge fault diameter   $\D^E_{a} =  \D_{(0,a)}$.
Let $G$ be a Cartesian graph bundle with fibre $F$ over the base graph $B$.
We show that  \\
(1)    $\D^V_{a+b+1}(G)\leq \D^V_{a}(F)+\D^V_{b}(B)$
when the graphs $F$ and $B$  are $k_F$-connected and $k_B$-connected,   $0< a < k_F$, $0< b < k_B$,
and  provided that $\D_{(a-1,1)}(F)\leq \D^{V}_{a} (F)$ and $\D_{(b-1,1)}(B)\leq \D^{V}_{b} (B)$ 
and \\
(2) $\D^E_{a+b+1}(G)\leq \D^E_{a}(F)+\D^E_{b}(B)$
when the graphs $F$ and $B$ are  $k_F$-edge connected and $k_B$-edge connected, $0\leq a < k_F$, $0\leq b < k_B$,
and provided  that $\D^E_{a}(F)\geq 2$ and $\D^E_{b}(B)\geq 2$.
\end{abstract}


\end{frontmatter}


\newtheorem{theorem}{Theorem}[section]
\newtheorem{corollary}[theorem]{Corollary}
\newtheorem{proposition}[theorem]{Proposition}
\newtheorem{lemma}[theorem]{Lemma}
\newtheorem{remark}[theorem]{Remark}
\newtheorem{definition}[theorem]{Definition}
\newtheorem{question}[theorem]{Question}
\newtheorem{example}[theorem]{Example}
\newtheorem{conjecture}[theorem]{Conjecture}


\section{Introduction}

The concept of fault diameter of Cartesian product graphs was first described in  \cite{1}, but the upper bound was wrong, 
as shown by  Xu, Xu and Hou who  provided a small counter example and  corrected the mistake \cite{0}.
More precisely,  denote by $\D^V_{a }(G)$ the fault diameter of a graph $G$, a maximum diameter of $G$ after deletion of any $a$ vertices,
and $G\Box H$ the Cartesian product of graphs $G$ and $H$.
 Xu, Xu and Hou proved \cite{0}
$$\D^V_{a+b+1}(G\Box H)\leq \D^V_{a}(G)+\D^V_{b}(H)+1
$$
while the claimed bound in \cite{1} was $ \D^V_{a}(G)+\D^V_{b}(H)$. 
(Our notation here slightly differs from notation used in  \cite{1} and \cite{0}.)
The result was later generalized to graph bundles in \cite{zer-ban} and generalized graph products (as defined by \cite{Bermond})  in \cite{KITAJCIdiameter}.
Here we 
show that in most cases of Cartesian graph bundles the bound can indeed be improved to the one claimed in \cite{1}.

Methods used involve the theory of mixed connectivity and  recent results on mixed fault diameters \cite{3fd,mpsv,fdsv,fdsv2}.
For completeness, we also give the analogous improved upper bound for edge fault diameter.

The rest of the paper is organized as follows.  In the next section we recall that the graph products and graph bundles often appear as 
practical interconnection network topologies because of some attractive properties they have.
In Section \ref{Preliminaries} we  provide general definitions, in particular of the connectivities.
Section  \ref{Bundles} introduces graph bundles and recalls relevant previous results.
The improved bounds are  proved in Section \ref{Results}. 

\section{Motivation - interconnection networks} 

 Graph products  and bundles belong to a class of frequently  studied interconnection network topologies.
For example meshes, tori, hypercubes and some of their generalizations are Cartesian products.
It is less known that  some other well-known interconnection network  topologies 
are  Cartesian graph bundles, for example  twisted hypercubes \cite{cull,efe} and 
multiplicative circulant graphs \cite{ivan}.  

In the design of large interconnection networks several factors have to be taken into account.
A usual constraint is that each  processor can be connected to a limited number of other processors and that
the delays in communication must not be too long. Furthermore, 
an interconnection network should be fault tolerant, because practical communication networks are exposed 
to failures of network components. 
Both failures of nodes and failures of connections between them happen and it is desirable that a network is robust in 
the sense that a limited number of failures does not break down the whole system.
A lot of work has been done on various aspects of network fault tolerance,
see for example the survey \cite{Bermond} and the more recent papers \cite{hung,sun,yin}.
In particular the fault diameter with faulty vertices, which was first studied in \cite{1}, and the edge fault diameter
have
been determined for many important networks recently \cite{zer-ban,zbedge,ZB,erves,day,4,7,0}. 
Usually either only edge faults or only vertex faults are considered, 
while the case when both edges and vertices may be faulty is studied rarely. 
For example, \cite{hung,sun} consider Hamiltonian properties assuming a combination of vertex and edge faults. 
In recent work on fault diameter of Cartesian graph products and bundles \cite{zer-ban,zbedge,ZB,erves}, 
analogous results were found for both fault diameter and edge fault diameter. 
However, the proofs for vertex and edge faults 
are independent, 
and our effort to see how results in one case may imply the others was  not successful.
A natural question is whether it is possible to design  a uniform theory
that covers simultaneous faults of vertices and edges.  
Some basic results on edge, vertex and mixed fault diameters for general graphs appear in \cite{3fd}.
In order to study the fault diameters of graph products and bundles  under mixed faults, it is important to 
understand generalized connectivities. 
Mixed connectivity which generalizes both vertex and edge connectivity, and some basic observations for any connected graph are given in \cite{mpsv}. 
%
 We are not aware of any earlier work on mixed connectivity. A closely related notion is  the connectivity pairs of a graph \cite{Harary}, 
but after Mader \cite{Mader} showed the claimed proof of generalized Menger's theorem is not valid, work on connectivity pairs seems to be very rare.

Upper
bounds for the mixed fault diameter of Cartesian graph bundles are  given in \cite{fdsv,fdsv2}
that  in some case also improve previously known results on vertex and edge fault diameters on these classes of 
Cartesian graph bundles \cite{zer-ban,erves}.  
However results in \cite{fdsv} address only the number of faults given by the connectivity of the fibre (plus one vertex), 
while the connectivity of the graph bundle can be much higher when the connectivity of the base graph is substantial, and
results in \cite{fdsv2} address only the number of faults given by the connectivity of the base graph (plus one vertex), 
while the connectivity of the graph bundle can be much higher when the connectivity of the fibre is substantial.
An upper bound for the mixed fault diameter that would take into account both types of faults remains to be an interesting open research problem.


\section{Preliminaries} \label{Preliminaries}
 
 A {\em simple graph}  $G=(V,E)$  is determined by a {\em vertex set} $V=V(G)$ 
and a set  $E=E(G)$ of (unordered) pairs of vertices, called  {\em edges}.
As usual, we  will use the short notation $uv$ for edge $\{u,v\}$.
For an edge $e=uv$ we call $u$ and $v$ its {\em endpoints}.
It is sometimes convenient to consider the union of  \emph{elements} of a graph, $S(G)= V(G) \cup E(G)$.
Given $X \subseteq S(G)$ then $S(G) \setminus X$ is a subset of elements of $G$.
However, note that in general  $S(G) \setminus X$ may not induce a graph. As we need notation for subgraphs with some missing (faulty) elements, 
we  formally  define $G\setminus X$, the subgraph of $G$ after deletion of $X$, as follows:
\begin{definition}
Let $X \subseteq S(G)$, and $X=X_E \cup X_V$, where $X_E \subseteq E(G)$ and $X_V \subseteq V(G)$.
Then $G\setminus X$ is the subgraph of $(V(G), E(G)\setminus X_E)$ induced on vertex set $V(G)\setminus X_V$. 
\end{definition}
A {\em walk} between vertices $x$ and $y$ is a sequence of vertices  and edges
$v_0,$ $e_1,$ $v_1,$ $e_2,$ $v_2,$ $\dots,$ $v_{k-1},$ $e_k,$ $v_k$
where $x=v_0$, $y=v_k$, and  $e_i = v_{i-1}v_i$ for each $i$.
A walk with all vertices distinct is called a {\em path},
and the vertices $v_0$ and $v_k$ are called the {\em endpoints} of the  path.
The {\em length} of a path $P$, denoted by $\ell (P)$, is the number of edges in $P$. 
The {\em distance} between vertices $x$ and $y$, denoted by $d_G(x,y)$, is the length of a shortest path 
between $x$ and $y$ in $G$.
If there is no path between $x$ and $y$ we write $d_G(x,y) = \infty$.
The {\em diameter} of a connected graph $G$, $\D(G)$, 
is the maximum distance between any two vertices in $G$. 
A path $P$ in $G$, defined by a sequence 
$x=v_0,e_1,v_1,e_2,v_2,\dots,v_{k-1},e_k,v_k=y$ can alternatively be seen as a subgraph of $G$
with $V(P) =\{v_0,v_1,v_2,\dots,v_k\}$ and $E(P) =\{e_1,e_2,\dots,e_k\}$.
Note that the reverse sequence gives rise to the same subgraph. 
Hence we use $P$ for a path either from $x$ to $y$ or from $y$ to $x$.
A graph is {\em connected} if there is a path between each pair of vertices, and is {\em disconnected} otherwise.
In particular, $K_1$ is by definition disconnected.
%
%
The  {\em connectivity} (or {\em vertex connectivity}) $\kappa (G)$ of a connected graph $G$, other than a complete graph, is the smallest number of vertices whose removal disconnects $G$.
For complete graphs is $\kappa (K_n)=n-1$. 
We say that $G$ is  {\em $k$-connected} (or {\em $k$-vertex connected}) for any $0<k \leq \kappa (G)$.
The  {\em edge connectivity} $\lambda (G)$ of a connected graph $G$, is the smallest number of edges whose removal disconnects $G$. 
A graph $G$ is said to be {\em $k$-edge connected}  for any $0<k \leq \lambda (G)$.
It is well known that
(see, for example, \cite{m}, page 224)
$\kappa (G) \leq \lambda (G) \leq \delta_G, $
where $\delta_G$ is smallest vertex degree of $G$. 
Thus if a graph $G$ is $k$-connected, then it is also $k$-edge connected. The reverse does not hold in general.

The mixed connectivity generalizes both vertex and edge connectivity \cite{mpsv,fdsv}.
Note that the definition used in \cite{fdsv} and  here slightly differs from the definition used in a previous work \cite{mpsv}.
%
{\begin{definition}
Let $G$ be any connected graph. A graph $G$ is
\emph{$(p,q)$+connected}, if $G$ remains connected after removal of any $p$ vertices and any $q$ edges. 
\end{definition}} 

We wish to remark that the mixed connectivity studied here is closely related to connectivity pairs as defined in \cite{Harary}.
Briefly speaking, a connectivity pair  of a graph is an ordered pair  $(k,\ell)$ of two integers such that there is some set of $k$ vertices 
and $\ell$ edges whose removal disconnects the graph and there is no set of $k-1$ vertices and $\ell$ edges  or of $k$ vertices and $\ell-1$ edges
with this property.
Clearly  $(k,\ell)$ is a connectivity pair of $G$ exactly when:
(1) $G$ is $(k-1,\ell)$+connected, 
(2) $G$ is $(k,\ell-1)$+connected,  and 
(3) $G$ is not  $(k,\ell)$+connected.
In fact, as shown in  \cite{mpsv}, (2) implies (1), so $(k,\ell)$ is a connectivity pair exactly when (2) and (3) hold.

From the definition we easily observe that any connected graph $G$ is $(0,0)$+ connected, 
$(p,0)$+connected for any $p<\kappa(G)$ and $(0,q)$+connected for any $q<\lambda(G)$.
In our notation $(i,0)$+connected is the same as $(i+1)$-connected, i.e. the graph remains connected 
after removal of any $i$ vertices.
Similarly, $(0,j)$+connected   means  $(j+1)$-edge connected, i.e. the graph remains con\-nected 
after removal of any $j$ edges.

Clearly, if $G$ is a $(p,q)$+connected graph, then $G$ is $(p^{\prime},q^{\prime})$+connected for any
$p^{\prime}\leq p$ and any $q^{\prime}\leq q$.
Furthermore, for any connected graph $G$ with $k<\kappa(G)$ faulty vertices, at least $k$ edges are not working.
Roughly speaking, graph $G$ remains connected if any faulty vertex in $G$ is replaced with a faulty edge. 
It is known \cite{mpsv}  
that if a graph $G$ is $(p,q)$+connected and $p>0$,  then $G$ is $(p-1,q+1)$+connected.
Hence for $p>0$ we have a chain of implications: 
$(p,q)$+connected $\Longrightarrow$ $(p-1,q+1)$+connected $\Longrightarrow\dots\Longrightarrow$ 
$(1,p+q-1)$+connected $\Longrightarrow$  $(0,p+q)$+connected, 
which 
generalizes the  well-known proposition 
that any $k$-connected graph is also $k$-edge connected. 
Therefore, a graph $G$ is $(p,q)$+connected if and only if $p<\kappa(G)$ and $p+q<\lambda(G)$. 

Note that by our definition the complete graph $K_n$, $n\geq 2$, is $(n-2,0)$+ connected, and hence $(i,j)$+connected for any $i+j\leq n-2$. Graph $K_2$ is $(0,0)$+connected, and mixed connectivity of $K_1$ is not defined.

If for a graph $G$ $\kappa(G)=\lambda(G) =k$, then $G$ is $(i,j)$+connected exactly when $i+j< k$.
However, if  $2\leq \kappa(G)<\lambda(G)$, the question whether $G$ is $(i,j)$+connected  for  $1\leq i<\kappa(G)\leq i+j < \lambda(G)$
is not trivial. 
The example below shows that in general the knowledge of 
$\kappa(G)$ and $\lambda(G)$ is not enough to decide whether  $G$ is  $(i,j)$+connected.

\begin{figure}[htb] 
\begin{center}
    \includegraphics[width=4.0in]{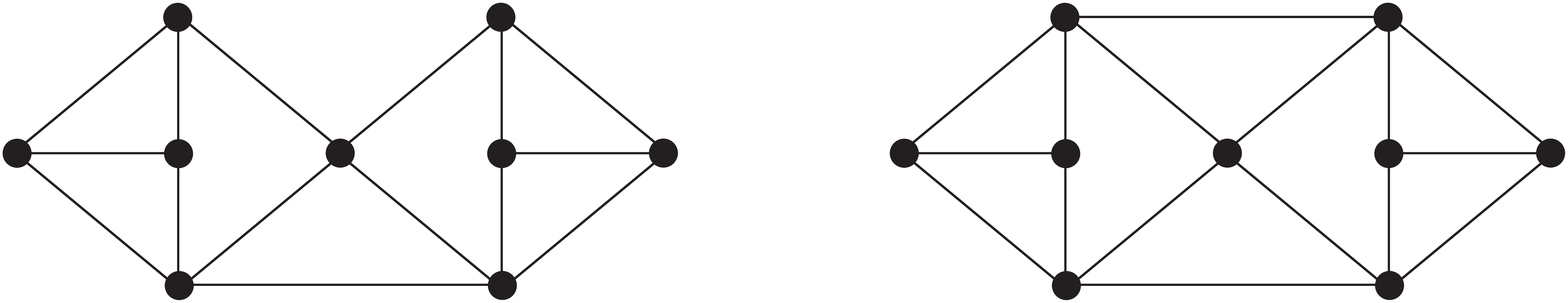}
    \caption{Graphs $G_1$ and $G_2$ from Example \ref{pr}.}
   \label{sl}
  \end{center}
\end{figure}

\begin{example} \label{pr}
For graphs on Fig. \ref{sl} we have $\kappa(G_1)=\kappa(G_2)=2$ and $\lambda(G_1)=\lambda(G_2)=3$. Both graphs are 
$(1,0)+$connected $\Longrightarrow$ $(0,1)+$connected, and $(0,2)+$ connected.
Graph $G_1$ is not $(1,1)+$connected, while graph $G_2$ is.
\end{example}

{\begin{definition}
Let $G$ be a $k$-edge connected graph and $0\leq a < k$. The {\em $a$-edge fault diameter} of $G$ is 
$$\D^E_a(G)=\max{\{\D(G\setminus X)\  | \  X\subseteq E(G),\left|X\right|=a\}}.$$ 
\end{definition}} 
 
{\begin{definition}
Let $G$ be a $k$-connected graph and $0\leq a < k$. The {\em $a$-fault diameter} (or {\em $a$-vertex fault diameter}) of $G$ is 
$$\D^V_a(G)=\max{\{\D(G\setminus X)\  | \  X\subseteq V(G),\left|X\right|=a\}}.$$ 
\end{definition}} 

Note that $\D^E_a(G)$ is the largest diameter among the diameters of subgraphs of $G$ with $a$ edges deleted, 
and $\D^V_a(G)$ is the largest diameter over all 
subgraphs of $G$ with $a$ vertices deleted.
In particular,  $\D^E_0(G)= \D^V_0(G)= \D(G)$, the diameter of $G$.
For  $p \geq \kappa (G)$
and for  $q \geq \lambda (G)$
we set $\D^V_p(G) = \infty$, $\D^E_q(G) = \infty$,  as some of the subgraphs are not vertex connected or edge connected, respectively.

It is known  \cite{3fd}
that for any connected graph $G$ the inequalities below hold.
\begin{enumerate}
\item $\D(G) =\D^E_0(G) \leq \D^E_1(G) \leq \D^E_2(G) \leq \ldots \leq \D^E_{\lambda (G)-1} (G)< \infty$. 
\item $\D(G) =\D^V_0(G) \leq \D^V_1(G) \leq \D^V_2(G) \leq \ldots \leq \D^V_{\kappa (G)-1} (G)< \infty$. 
\end{enumerate}

\begin{definition} \label{MFD}
Let $G$ be a $(p,q)$+connected graph.
The {\em $(p,q)$-mixed fault dia\-meter} of $G$ is 
$$\D_{(p,q)} (G)=\max{\{\D(G\setminus (X\cup Y))\  | \ X\subseteq V(G), Y\subseteq E(G), |X|=p, |Y|=q\}}.$$ 
\end{definition}

Note that by Definition \ref{MFD} the endpoints of edges of set $Y$ can be in $X$. In this case we may get 
the same subgraph of $G$ by deleting $p$  vertices and fewer than $q$ edges. 
It is however not difficult to see that the diameter of such subgraph is smaller than or equal to the 
diameter of some subgraph of $G$ where exactly $p$ vertices and exactly $q$ edges are deleted. 
So the condition that the endpoints of edges of set $Y$ are not in $X$ need not to be included in Definition \ref{MFD}.
The mixed fault diameter $\D_{(p,q)} (G)$ is the largest diameter among the diameters of all subgraphs obtained from $G$ by deleting $p$ vertices and $q$ edges,
hence $\D_{(0,0)}(G)= \D(G)$, $\D_{(0,a)}(G)=\D^E_a(G)$ and $\D_{(a,0)}(G)=\D^V_a(G)$.

Let $\H_a^V=\{G\setminus X\  | \  X\subseteq V(G), |X|=a\}$
and
$\H_b^E=\{G\setminus X\  | \  X\subseteq E(G), |X|=b\}.$
It is easy to see that
\begin{enumerate}
\item $\max{\{\D^E_b(H) \ | \ H\in \H_a^V\}} = \D_{(a,b)}(G)$,
\item $\max{\{\D^V_a(H) \ | \ H\in \H_b^E\}} = \D_{(a,b)}(G)$.
\end{enumerate}

In previous work \cite{3fd} on vertex, edge and mixed fault diameters of connected graphs the following theorem has been proved. 

{\begin{theorem}\label{main_mixed} 
Let $G$ be $\left(  p,q\right)$+connected graph and $p>0$. 
\begin{itemize}
\item If $q>0$, then $\D^{E}_{p+q} (G) \leq \D_{(1,p+q-1)}(G) \leq \dots \leq \D_{(p,q)} (G)$. 
\item If $q=0$, then $\D^{E}_{p} (G) \leq \D_{(1,p-1)} (G)\leq \dots \leq \D_{(p-1,1)} (G)\leq \D^{V}_{p} (G)+1$.
\end{itemize}
\end{theorem}}

Note that for $(p+1)$-connected graph $G$, $p>0$, we have     either 
$$\D^{E}_{p} (G) \leq \D_{(1,p-1)} (G)\leq \dots \leq \D_{(p-1,1)} (G)\leq \D^{V}_{p} (G)$$ or 
$$\D^{E}_{p} (G) \leq \D_{(1,p-1)} (G)\leq \dots \leq \D_{(p-1,1)} (G)= \D^{V}_{p} (G)+1.$$

For example, complete graphs, complete bipartite graphs, and  cycles are graphs with 
$\D_{(p-1,1)} (G)= \D^{V}_{p}(G)+1$ for all meaningful of values of $p$.
More examples of both types of graphs can be found in \cite{3fd}.

 
\section{Fault diameters of Cartesian graph bundles}
\label{Bundles} 

Cartesian graph bundles are a generalization of Cartesian graph products, first studied in \cite{PiVr,ST_P_V}. 
Let $G_1$ and $G_2$ be graphs. The {\em Cartesian product} of graphs $G_1$ and $G_2$, $G=G_1\Box G_2$,
is defined on the vertex set $V(G_1)\times V(G_2)$.
Vertices $(u_1,v_1)$ and $(u_2,v_2)$ are 
adjacent if either $u_1u_2\in E(G_1)$ and $v_1=v_2$ or $v_1v_2\in E(G_2)$ and $u_1=u_2$. 
For further reading on graph products we recommend  \cite{3}. 

{\begin{definition}
Let $B$ and $F$ be graphs. A graph $G$ is a {\em Cartesian graph bundle with fibre $F$ over the base graph} $B$ if 
there is a {\em graph map} $p:G\rightarrow B$ such that for each vertex $v\in V(B)$, $p^{-1}(\{v\})$ is isomorphic to $F$, 
and for each edge $e=uv\in E(B)$, $p^{-1}(\{e\})$ is isomorphic to $F \Box K_2$.
\end{definition}} 
More precisely, the mapping $p:G\rightarrow B$  maps graph elements of $G$ to graph elements of $B$, 
i.e. $p: V(G) \cup E(G) \rightarrow  V(B) \cup E(B)$. In particular, here we also assume that 
the vertices of $G$ are mapped to vertices of $B$ and the edges of $G$ are mapped either to vertices or to edges of $B$. 
We say an edge $e\in E(G)$ is {\em degenerate} if $p(e)$ is a vertex. Otherwise we call it {\em nondegenerate}.
The mapping $p$ will also be called the {\em projection} (of the bundle $G$ to its base $B$).
Note that  each edge $e=uv \in E(B)$ naturally induces an isomorphism $\varphi_e :p^{-1}(\{u\})\rightarrow p^{-1}(\{v\})$
between two fibres. 
It may be interesting to note that while it is well-known that a graph can have only one representation  as a product
(up to isomorphism and up to the order of factors) \cite{3},
there may be many different graph bundle representations of the same graph \cite{ZmZe2002b}.
Here we assume that the bundle representation is given.
Note that in some cases finding a representation of $G$ as a graph bundle can be found in polynomial time 
\cite{ImPiZe,ZmZe2000,ZmZe2001a,ZmZe2002b,ZmZe2002a,directed}. 
For example, one of the easy classes are the Cartesian graph bundles over triangle-free base \cite{ImPiZe}.
Note that a graph bundle over a tree $T$ (as a base graph) with fibre $F$ is isomorphic to the 
Cartesian product $T\Box F$ (not difficult to see, appears already in  \cite{PiVr}), 
i.e. we can assume that all isomorphisms $\varphi_e$ are identities.
For a later reference note that for any path $P\subseteq B$,
$p^{-1}(P)$ is a  Cartesian graph bundle over the path $P$, and  
one can  define coordinates in the product $P\Box F$ in a natural way.
 
%
%
%
%
%
%
%
In recent work on fault diameter of Cartesian graph products and bundles \cite{zer-ban,zbedge,ZB,erves}, 
analogous results were found for both fault diameter and edge fault diameter. 

\begin{theorem} \cite{zer-ban}
\label{VPsv}
Let $F$ and $B$ be $k_F$-connected and $k_B$-connected graphs, respectively, $0\leq a < k_F$,  
$0\leq b < k_B$, and $G$ a Cartesian bundle with fibre $F$ over the base graph $B$. Then 
$$\D^V_{a+b+1}(G)\leq \D^V_{a}(F)+\D^V_{b}(B)+1.$$
\end{theorem} 

\begin{theorem} \cite{erves}
\label{EPsv}
Let $F$ and $B$ be $k_F$-edge connected and $k_B$-edge connected graphs, respectively, $0\leq a < k_F$, $0\leq b < k_B$, and $G$ a Cartesian bundle with fibre $F$ over the base graph $B$. Then 
$$\D^E_{a+b+1}(G)\leq \D^E_{a}(F)+\D^E_{b}(B)+1.$$
\end{theorem} 

Before writing a theorem on bounds for the mixed fault diameter we recall a theorem on mixed connectivity. 

\begin{theorem} \cite{mpsv} 
\label{mpsv2} Let $G$ be a Cartesian graph bundle with fibre $F$ over the base graph $B$, graph $F$ be $(p_{F},q_{F})$+connected and
graph $B$ be $(p_{B},q_{B})$+connec\-ted.
Then Cartesian graph bundle $G$ is $(p_{F}+p_{B}+1,q_{F}+q_{B})$+connected.
\end{theorem}

In recent work \cite{fdsv,fdsv2}, an upper bound for the mixed fault diameter of Cartesian graph bundles, $\D_{(p+1,q)} (G)$, in terms of mixed fault diameter of the fibre and diameter of the base graph and in terms of diameter of the fibre and mixed fault diameter of the base graph, respectively, 
is given.

\begin{theorem}  \cite{fdsv} 
\label{MPsvBsp}
Let $G$ be a Cartesian graph bundle with fibre $F$ over the base graph $B$, 
where graph $F$ is $(p,q)$+connected, $p+q>0$, 
and $B$ is a connected graph with diameter $\D(B)>1.$
Then we have:
\begin{itemize}
\item If $q>0$, then $\D_{(p+1,q)} (G)\leq \D_{(p,q)}(F)+ \D(B).$
\item If $q=0$, then $\D^{V}_{p+1} (G)\leq \max \{ \D^{V}_{p} (F), \D_{(p-1,1)}(F)\}+ \D(B).$
\end{itemize}
\end{theorem} 

Theorem \ref{MPsvBsp} improves results \ref{VPsv} and  \ref{EPsv} for $a>0$ and $b=0$. \\
Let $G$ be a Cartesian graph bundle with fibre $F$ over the connected base graph $B$ with  diameter $\D(B)>1$, and let $a>0$.
If graph $F$ is $(a+1)$-connected, i.e. $(a,0)$+connected, then  
by theorem \ref{MPsvBsp} we have an upper bound for the vertex fault diameter  
$\D^{V}_{a+1} (G) \leq \D^{V}_{a} (F)+ \D(B)+1$ for any graph $F$.
Similarly, 
$\D^{V}_{a+1} (G)\leq  \D^{V}_{a} (F)+ \D(B) $ if  $\D_{(a-1,1)}(F) \leq \D^{V}_{a} (F)$ holds.\\
If graph $F$ is $(a+1)$-edge connected, i.e. $(0,a)$+connected, then
by theorems \ref{main_mixed} and \ref{MPsvBsp} we have an upper bound for the edge fault diameter  
$\D^{E}_{a+1} (G) \leq \D_{(1,a)}(G)\leq   \D^{E}_{a} (F)+ \D(B)$.

\begin{theorem}  \cite{fdsv2} 
\label{MPsvFsp}
Let $G$ be a Cartesian graph bundle with fibre $F$ over the base graph $B$,  graph $F$ be a connected graph with diameter $\D(F)>1$, and graph $B$ be   $\left(p,q\right)$+connected, $p+q>0$. 
Then we have:

\begin{itemize}
\item If $q>0$, then $\D_{(p+1,q)} (G)\leq \D(F) + \D_{(p,q)} (B).$ 
\item If $q=0$, then $\D^{V}_{p+1} (G)\leq \D(F) + \max \{ \D^{V}_{p} (B), \D_{(p-1,1)}(B)\}.$
\end{itemize}
\end{theorem}

Theorem \ref{MPsvFsp} improves results \ref{VPsv} and  \ref{EPsv} for $a=0$  and $b>0$. \\
Let $G$ be a Cartesian graph bundle with fibre $F$ over the base graph $B$,  graph $F$ be a connected graph with diameter $\D(F)>1$, and let $b>0$.
If graph $B$ is $(b+1)$-connected, i.e. $(b,0)$+connected, then  
by  Theorem \ref{MPsvFsp} we have an upper bound for the vertex fault diameter 
$\D^{V}_{b+1} (G)\leq \D(F) + \D^{V}_{b} (B)+1$ 
 for any graph $B$.
Similarly, 
$\D^{V}_{b+1} (G)\leq \D(F) + \D^{V}_{b} (B)$
if  $\D_{(b-1,1)}(B) \leq \D^{V}_{b} (B)$ holds. \\
If graph $B$  is $(b+1)$-edge connected, i.e. $(0,b)$+connected, then
by  theorems \ref{main_mixed} and \ref{MPsvBsp} we have an upper bound for the edge fault diameter  
$\D^{E}_{b+1} (G) \leq \D_{(1,b)}(G)\leq  \D (F)+ \D^{E}_{b} (B)$.

In the case when $a=b=0$ the fault diameter is determined exactly.
\begin{proposition} 
\cite{fdsv}
Let $G$ be a Cartesian graph bundle with fibre $F$ over the base graph $B$, and graphs $F$ and $B$ be connected graphs with diameters $\D(F)>1$ and $\D(B)>1$.
Then $$\D^{V}_{1} (G)= \D^{E}_{1} (G)=\D(G)=\D(F) + \D(B).$$
\end{proposition}

In other words, the diameter of a nontrivial Cartesian graph  bundle does not change when one element is faulty.

Here we improve  results  of theorems \ref{VPsv} and  \ref{EPsv} for positive $a$ and $b$.

\section{The results - improved bounds} \label{Results}

Before stating and proving the main theorems, 
we introduce some notation used in this section. \\
Let $G$ be a Cartesian graph bundle with fibre $F$ over the base graph $B$.
The {\em fibre of vertex} $x \in V(G)$ is denoted by $F_x$, formally, $F_x = p^{-1}(\{p(x)\})$. We will also use notation $F(u)$ 
for the fibre of the vertex $u \in V(B)$, i.e. $F(u) = p^{-1}(\{u\})$. Note that $F_x = F(p(x))$. 
We  will  also  use  shorter  notation $x\in F(u)$ for  $x\in V(F(u))$. \\ 
Let $u,v\in V(B)$ be distinct vertices, and
$Q$ be a path from $u$ to $v$ in $B$, and $x\in F(u)$. Then the {\em lift of the path $Q$ to the vertex} $x\in V(G)$, 
$\tilde Q_x$, is the path from $x\in F(u)$ to a vertex in $F(v)$, 
such that $p(\tilde Q_x)=Q$ and $\ell(\tilde Q_x)=\ell(Q)$.
Let $x,x'\in F(u)$. Then $\tilde Q_x$ and $\tilde Q_{x'}$ have different endpoints in $F(v)$ and are disjoint paths if 
and only if $x\neq x'$.
In fact,  two lifts  $\tilde Q_x$ and $\tilde Q_{x'}$ are either disjoint  $\tilde Q_x \cap \tilde Q_{x'} = \emptyset$
or equal,  $\tilde Q_x = \tilde Q_{x'}$.
We will also use notation $\tilde Q$ for lifts of path $Q$ to any vertex in $F(u)$. \\
Let $Q$ be a path from $u$ to $v$ and $e=uw\in E(Q)$. We will use notation $Q \setminus e$ for the subpath from $w$ to $v$,  i.e.  $Q \setminus e = Q \setminus \{u,e\}= Q \setminus \{u\}$.\\
Let $G$ be a graph and $X\subseteq S(G)$ be a set of elements of $G$. A path $P$ from a vertex $x$ to a vertex $y$ 
{\em avoids} $X$ in $G$, 
if $S(P)\cap X=\emptyset $, and it {\em internally avoids} $X$, 
if $(S(P)\setminus \{x,y\})\cap X=\emptyset$.

\subsection{Vertex fault diameter of Cartesian graph bundles}

\begin{theorem} 
\label{VFD} 
Let $G$ be a Cartesian graph bundle with fibre $F$ over the base graph $B$, graphs $F$ and $B$ be $k_F$-connected and $k_B$-connected, respectively, and let $0< a < k_F$, $0< b < k_B$. 
If for fault diameters of graphs $F$ and $B$, $\D_{(a-1,1)}(F)\leq \D^{V}_{a} (F)$ and $\D_{(b-1,1)}(B)\leq \D^{V}_{b} (B)$ hold
then 
$$\D^V_{a+b+1}(G)\leq \D^V_{a}(F)+\D^V_{b}(B).$$
\end{theorem}

\begin{pf*}{Proof.} 
Let $G$ be a Cartesian graph bundle with fibre $F$ over the base graph $B$, graph $F$  be $(a+1)$-connected, $ a >0$, graph $B$ be $(b+1)$-connected, $b >0$, and let $\D_{(a-1,1)}(F)\leq \D^{V}_{a} (F)$, $\D_{(b-1,1)}(B)\leq \D^{V}_{b} (B)$.   
Then $\D^{V}_{a} (F) \geq 2$, $\D^{V}_{b} (B) \geq 2$, and Cartesian bundle $G$ is  $(a+b+2)$-co\-nnec\-ted.
Let $X\subseteq V(G)$ be a set of faulty vertices, $\left|X\right|= a+b+1$, 
and let $x,y\in V(G)\setminus X$ be two distinct nonfaulty vertices in $G$. 
We shall consider the distance $d_{G\setminus X}(x,y)$.
\begin{itemize}
\item	Suppose first that $x$ and $y$ are in the same fibre, i.e. $p(x) = p(y)$.\\
If $\left| X \cap V(F_x)\right|\leq a$, then $d_{G\setminus X}(x,y)\leq \D^{V}_{a} (F)$.\\
If $\left| X \cap V(F_x)\right| > a$, then outside of fibre $F_x$ there are at most $b$ 
faulty vertices.
As graph $B$ is $(b+1)$-connected, there are at least $b+1$ neighbors of vertex $p(x)$ in $B$. 
Therefore there exist a neighbor $v$ of vertex $p(x)$ in $B$,   such that $\left| X \cap F(v)\right|=0$, and there is a path $x\rightarrow x^{\prime} \stackrel{P}{\rightarrow}y^{\prime}\rightarrow y$, which
avoids $X$, where $x^{\prime},y^{\prime} \in F(v)$ and $\ell(P)\leq \D(F)$. 
Thus $d_{G\setminus X}(x,y)\leq 1+\D (F)+1\leq  \D^{V}_{a} (F)+\D^{V}_{b} (B)$. 

\item	Now assume that $x$ and $y$  are in distinct fibres, i.e. $p(x) \neq p(y)$.
Let $X_B=\{v \in V(B)\setminus \{p(x),p(y)\}  ;  \left| X \cap F(v)\right|>0\}$.
We distinguish two cases.

\begin{enumerate}
\item
If $\left|X_B\right| \geq b$, then let $X_B^{\prime} \subseteq X_B$
 be an arbitrary subset of $X_B$ with 
 $\left| X_B^{\prime} \right|=b$.
The subgraph $B \setminus X_B^{\prime}$ is  a connected graph and there exists a path $Q$  in  $B \setminus X_B^{\prime}$ 
from $p(x)$ to $p(y)$ with $\ell(Q)\leq \D^{V}_{b} (B)$.
 In $p^{-1}(Q)=F\Box Q$ there are at most $a+1$ faulty vertices.
Let $x^{\prime} \in F_y$ be the endpoint of the path $\tilde{Q}_x$, the lift of $Q$. 
We distinguish two cases.

\begin{enumerate}
\item
If $x^{\prime}=y$, then $\tilde{Q}_x$ is a path from $x$ to $y$ in $G$.  
If $\tilde{Q}_x$ avoids $X$, then $d_{G\setminus X} (x,y)\leq \ell(Q)\leq \D^{V}_{b} (B)$.
If $\tilde{Q}_x$ does not avoid $X$,
then there are at most $a$ faulty vertices outside of the path $\tilde{Q}_x$ in $F\Box Q$.
As the  graph $F$ is $(a+1)$-connected, there are at least $a+1$ neighbors of $x$ in $F_x$.
Since there are more neighbors than faulty vertices (outside of $\tilde{Q}_x$ in $F\Box Q$), there exists 
a neighbor $v\in V(F_x)$ of $x$, such that the lift $\tilde{Q}_v$ avoids $X$.
The endpoint of the path $\tilde{Q}_v$ in fibre $F_y$ is a neighbor of $y$, therefore 
$d_{G\setminus X}(x,y) \leq 1+ \ell(Q) +1\leq \D^{V}_{a} (F)+ \D^{V}_{b} (B)$.
 
\item Let $x^{\prime}\neq y$. 
If $\left|V(F_x) \cap X \right|=a+1$ or $\left|V(F_y) \cap X \right|=a+1$, then obviously  $d_{G\setminus X}(x,y) \leq \ell(Q) + \D(F)\leq \D^{V}_{b} (B)+ \D^{V}_{a} (F)$.\\
Now assume $\left|V(F_x) \cap X \right|\leq a$ and $\left|V(F_y) \cap X \right|\leq a$.
If $\tilde{Q}_x$ or $\tilde{Q}_y$ avoids $X$, then $d_{G\setminus X}(x,y) \leq \ell(Q) + \D^{V}_{a} (F)\leq  \D^{V}_{b} (B)+\D^{V}_{a} (F)$.
Suppose that paths $\tilde{Q}_x$ and $\tilde{Q}_y$ do not avoid $X$. Then there are at most $a-1$ faulty vertices outside of paths $\tilde{Q}_x$ and $\tilde{Q}_y$ in $F \Box Q$.
Let $X^{\prime}  \subseteq V(F_y)$ be defined as $X^{\prime} =\{v\in V(F_y)\setminus\{x^{\prime}, y\}, \left|\tilde{Q}_v \cap X \right|>0 \}$.
Then $ \left|X^{\prime} \right| \leq a-1$.
There is a path $P$ from $x^{\prime}$ to $y$ in $F_y\setminus X^{\prime}$ of length $\ell(P)\leq \D^V_{a-1}(F)\leq \D^V_{a}(F)$. 
Note that the path $P$ internally avoids $X$.
If $x^{\prime}$ and $y$ are not adjacent, then $\ell(P)\geq 2$. 
For the  neighbor $v^{\prime}$ of $x^{\prime}$ on the path $P$, $e^{\prime}=x^{\prime}v^{\prime}\subset P$,
the lift $\tilde{Q}_{v^{\prime}}$ avoids $X$. Let $v\in V(F_x)$ be the endpoint of the lift $\tilde{Q}_{v^{\prime}}$. 
Then the path  $x\rightarrow v\stackrel{\tilde{Q}}{\rightarrow}v^{\prime}\stackrel{P\setminus e^{\prime}}{\rightarrow}y$ avoids $X$, therefore
$d_{G\setminus X}(x,y) \leq 1+\ell(Q) + \D^{V}_{a} (F)-1 \leq \D^{V}_{a} (F)+\D^{V}_{b} (B)$.\\
If $x^{\prime}$ and $y$ are adjacent, then remove from $F_y$ the set of vertices $X^{\prime}$ and the edge $e=x^{\prime}y$.
There is a path $P^{\prime}$ from $x^{\prime}$ to $y$ in $F_y\setminus (X^{\prime} \cup \{e\})$ of  length  
$2\leq \ell(P^{\prime})\leq \D_{(a-1,1)}(F)$, that internally avoids $X$. 
As before, for the neighbor $w^{\prime}$ of $x^{\prime}$ on the path $P^{\prime}$ the lift $\tilde{Q}_{w^{\prime}}$ avoids $X$. Therefore
$d_{G\setminus X}(x,y) \leq 1+\ell(Q) + \D_{(a-1,1)}(F)-1 \leq   \D^{V}_{a} (F)+\D^{V}_{b} (B)$.
\end{enumerate}

\item If $\left|X_B\right| < b$, then the subgraph $B \setminus X_B$ is (at least) $2$-connected, thus also $2$-edge connected.
If  the  vertex $p(y)$ is not a neighbor of $p(x)$, then there is a path 
$Q$ from $p(x)$ to $p(y)$ in $B$ with $2\leq \ell(Q)\leq \D^{V}_{b-1} (B)\leq \D^{V}_{b} (B)$ that internally avoids $X_B$. 
Let $v\in V(Q)$ be a neighbor of $p(x)$, $e^{\prime}=p(x)v$.
Then there is a path $x\rightarrow x^{\prime} \stackrel{P}{\rightarrow}y^{\prime}\stackrel{\tilde{Q}\setminus e^{\prime}}{\rightarrow} y$, which
avoids $X$, where  $x^{\prime},y^{\prime} \in F(v)$ and $\ell(P)\leq \D(F)$. 
Thus $d_{G\setminus X}(x,y)\leq 1+\D (F)+\D^{V}_{b} (B)-1 \leq \D^{V}_{a}(F)+\D^{V}_{b} (B)$. \\
If $e=p(x)p(y) \in E(B)$, then $B \setminus (X_B \cup \{e\})$ is a connected graph and there is a path $Q^{\prime}$ from $p(x)$ to $p(y)$ with $2\leq \ell(Q^{\prime})\leq \D_{(b-1,1)} (B)$ that internally avoids $X_B$.
Similarly as before we have $d_{G\setminus X}(x,y)\leq 1+\D (F)+\D_{(b-1,1)} (B)-1  \leq \D^{V}_{a}(F)+\D^{V}_{b} (B)$. \qed
\end{enumerate}
\end{itemize}
 \end{pf*}

Theorem \ref{VFD} improves Theorem \ref{VPsv} on  the class of Cartesian graph bundles 
for which both, the fiber and the base graph, are at least 2-connected. 
Theorem \ref{VFD} also improves result of \cite{0} on the  Cartesian graph products with at least 2-connected factors.
The next example shows that the bound of Theorem \ref{VFD} is tight.

\begin{example} 
Let $F=B=K_4\setminus \{e\}$.
Then graph  $F$ is $2$-connected and
$\D^{E}_{1}(F)=\D^{V}_{1}(F) =2$.
The vertex fault diameter of Cartesian graph product $F\Box F$ on Fig. \ref{K4-exK4-e} is  
$\D^{V}_{3}(F\Box F)= \D^{V}_{1} (F)+ \D^{V}_{1} (F)=4.$
%
\begin{figure}[hhhtb]
 \begin{center}
    \includegraphics[width=1.0in]{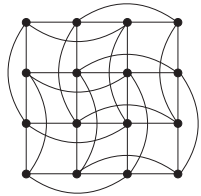}
    \caption{Cartesian graph product of two factors $K_4\setminus \{e\}$.}
    \label{K4-exK4-e}
  \end{center}
\end{figure}
\end{example} 
\begin{example} \label{V2}
Cycle $C_{4}$ is $2$-connected graph and
$\D^{E}_{1}(C_{4})=\D^{V}_{1}(C_{4})+1 =3$.
The vertex fault diameter of Cartesian graph bundle $G$ with fibre $C_4$ over base graph $C_4$ on Fig. \ref{iliac} is  
$\D^{V}_{3}(G)= \D^{V}_{1} (C_4)+ \D^{V}_{1} (C_4)+1=5.$
%
\begin{figure}[h!]
 \begin{center}
    \includegraphics[width=5.0in]{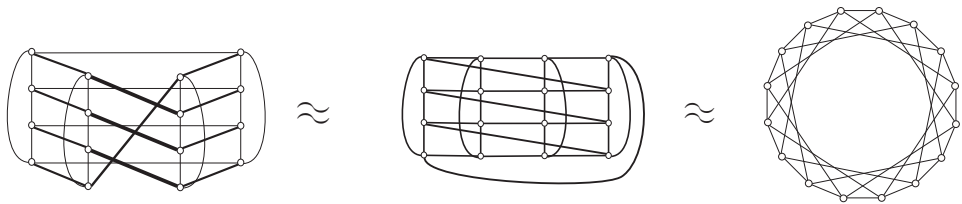}
    \caption{Twisted torus: Cartesian graph bundle with fibre $C_4$ over base $C_4$.}
    \label{iliac}
  \end{center}
\end{figure}
\end{example} 
 
It is less known that graph bundles also appear as computer topologies. 
A well known example is the twisted torus on Fig. \ref{iliac}.
Cartesian graph bundle with fibre $C_4$ over base $C_4$ is the ILLIAC IV architecture \cite{ILLIAC},
a famous supercomputer that inspired some modern multicomputer architectures.
It may be interesting to note that the original design was a graph bundle with fibre $C_8$ over base $C_8$, 
but due to high cost a smaller version was build \cite{computermuseum}.


\subsection{Edge fault diameter of Cartesian graph bundles}

Let $G$ be a $k$-edge connected graph and $0\leq a < k$. 
Note  that if $a>0$ then $\D^E_{a}(G)\geq 2$ for any graph $G$.
More precisely, $\D^E_{a}(G)\geq 2$ if $a>0$ or ($a=0$ and $G$ is not a complete graph).
Furthermore, $\D^E_{a}(G)=1$ if and only if $a=0$ and  $G$ is a complete graph.

\begin{theorem} 
\label{EFD} 
Let $G$ be a Cartesian graph bundle with fibre $F$ over the base graph $B$, graphs $F$ and $B$ be $k_F$-edge connected and $k_B$-edge connected, respectively, and let $0\leq a < k_F$, $0\leq b < k_B$. If for edge fault diameters of graphs $F$ and $B$, $\D^E_{a}(F)\geq 2$ and $\D^E_{b}(B)\geq 2$ hold then 
$$\D^E_{a+b+1}(G)\leq \D^E_{a}(F)+\D^E_{b}(B).$$
\end{theorem}

\begin{pf*}{Proof.} 
Let $G$ be a Cartesian graph bundle with fibre $F$ over the base graph $B$, the graph $F$  be $(a+1)$-edge connected, $\D^{E}_{a} (F) \geq 2$, 
and the graph $B$ be $(b+1)$-edge connected, $\D^{E}_{b} (B) \geq 2$.   
Then the  Cartesian bundle $G$ is  $(a+b+2)$-edge co\-nnec\-ted.
Let $Y\subseteq E(G)$ be the set of faulty edges, $\left|Y\right|=  a+b+1$. 
Denote the set of degenerate edges in $Y$ by $Y_D$, and the set of nondegenerate edges by $Y_N$, 
$Y=Y_N\cup Y_D$, $p(Y_D)\subseteq V(B)$, $p(Y_N)\subseteq E(B)$.  
Let $x,y\in V(G)$ be two arbitrary distinct vertices in $G$. 
We shall consider the distance $d_{G\setminus Y}(x,y)$.

\begin{itemize}
\item	Suppose first that $x$ and $y$ are in the same fibre, i.e. $p(x) = p(y)$.\\
If $\left| Y_D \cap E(F_x)\right|\leq a$, then $d_{G\setminus Y}(x,y)\leq \D^{E}_{a} (F)$.\\
If $\left| Y_D \cap E(F_x)\right| > a$, then outside of the fibre $F_x$ there are at most $b$ 
faulty edges.
As the graph $B$ is $(b+1)$-edge connected, there are at least $b+1$ neighbors of the  vertex $p(x)$ in $B$. 
Therefore there exist a neighbor $v$ of vertex $p(x)$ in $B$, $e=p(x)v \in E(B)$,  such that $\left| Y_D \cap F(v)\right|=0$ and $p(e) \notin p(Y_N)$, and
hence  there is a path $x\rightarrow x^{\prime} \stackrel{P}{\rightarrow}y^{\prime}\rightarrow y$, which
avoids $Y$, where $x^{\prime},y^{\prime} \in F(v)$ and $\ell(P)\leq \D(F)$. 
Thus $d_{G\setminus Y}(x,y)\leq 1+\D (F)+1\leq  \D^{E}_{a} (F)+\D^{E}_{b} (B)$.

\item	Now assume that $x$ and $y$  are in distinct fibres, i.e. $p(x) \neq p(y)$.
We distinguish two cases.

\begin{enumerate}
\item
If $\left|Y_N\right| \geq b$, then let $Y_N^{\prime} \subseteq Y_N$
 be an arbitrary subset of $Y_N$ with 
 $\left| Y_N^{\prime} \right|=b$.
The subgraph $B \setminus p(Y_N^{\prime})$ is a connected graph and there exists a path $Q$ from $p(x)$ to $p(y)$ with $\ell(Q)\leq \D^{E}_{b} (B)$.
 In $p^{-1}(Q)=F\Box Q$ there are at most $a+1$ faulty edges.
Let $x^{\prime} \in F_y$ be the endpoint of the path $\tilde{Q}_x$, the lift of $Q$. 
We distinguish two cases.

\begin{enumerate}
\item
If $x^{\prime}=y$, then $\tilde{Q}_x$ is a path from $x$ to $y$ in $G$.  
If $\tilde{Q}_x$ avoids $Y$, then $d_{G\setminus Y} (x,y)\leq \ell(Q)\leq \D^{E}_{b} (B)$.
If $\tilde{Q}_x$ does not avoid $Y$,
then there are at most $a$ faulty edges outside of the path $\tilde{Q}_x$ in $F\Box Q$.
As the graph $F$ is $(a+1)$-edge connected, there are at least $a+1$ neighbors of $x$ in $F_x$.
Since there are more neighbors than faulty edges (outside of $\tilde{Q}_x$ in $F\Box Q$) there exist  
a neighbor $s\in V(F_x)$ of $x$, such that the path $x\rightarrow s \stackrel{\tilde{Q}}{\rightarrow}s^{\prime}\rightarrow y$ avoids $Y$,
where $s^{\prime} \in V(F_y)$ is a neighbor of $y$, therefore 
$d_{G\setminus X}(x,y) \leq 1+ \ell(Q) +1\leq \D^{E}_{a} (F)+ \D^{E}_{b} (B)$.
 
\item Let $x^{\prime}\neq y$. 
If $\left|E(F_x) \cap Y \right|=a+1$ or $\left|E(F_y) \cap Y \right|=a+1$, then obviously  $d_{G\setminus Y}(x,y) \leq \ell(Q) + \D(F)\leq \D^{E}_{b} (B)+ \D^{E}_{a} (F)$.\\
Now let $\left|E(F_x) \cap Y \right|\leq a$ and $\left|E(F_y) \cap Y \right|\leq a$.
If $\tilde{Q}_x$ or $\tilde{Q}_y$ avoid $Y$, then $d_{G\setminus Y}(x,y) \leq \ell(Q) + \D^{E}_{a} (F)\leq  \D^{E}_{b} (B)+\D^{E}_{a} (F)$.
Suppose that paths $\tilde{Q}_x$ and $\tilde{Q}_y$ do not avoid $Y$. Then there are at most $a-1$ faulty edges outside of paths $\tilde{Q}_x$ and $\tilde{Q}_y$ in $F \Box Q$.
Let $Y_D^{\prime}  \subseteq E(F_y)$ be set of edges from $x^{\prime}$ to such neighbors $v^{\prime}_i\in V(F_y)$ for which the paths 
$v^{\prime}_i \stackrel{\tilde{Q}}{\rightarrow}v_i\rightarrow x$ do not avoid faulty edges,
$Y_D^{\prime} =\{e=x^{\prime}v^{\prime}\in E(F_y); \left|(\tilde{Q}_v^{\prime} \cup vx) \cap Y \right|>0, v= \tilde{Q}_v^{\prime} \cap F_x\}$.
Note that if $x^{\prime}$ is a neighbor of $y$ then  $x^{\prime}y \in Y_D^{\prime}$.
Then the subgraph $F_y \setminus (Y_D^{\prime} \cup Y_D)$ is a connected graph as there are at most $a+1$ faulty edges in $p^{-1}(Q)=F\Box Q$ and
 $\tilde{Q}_x$ does not avoid $Y$.
Therefore there is a path $P$ from $x^{\prime}$ to $y$ in $F_y\setminus (Y_D^{\prime} \cup Y_D)$ of length $2\leq \ell(P)\leq \D^E_{a}(F)$, which avoids $Y$ and for 
the  neighbor $v^{\prime}$ of $x^{\prime}$ on the path $P$, the lift $\tilde{Q}_{v^{\prime}}$ avoids $Y$. Let  $v= \tilde{Q}_{v^{\prime}} \cap F_x$. Then $vx \notin Y$, and
the path  $x\rightarrow v\stackrel{\tilde{Q}}{\rightarrow}v^{\prime}\stackrel{P\setminus e^{\prime}}{\rightarrow}y$ avoids $Y$, thus
$d_{G\setminus Y}(x,y) \leq 1+\ell(Q) + \ell(P)-1 \leq \D^{E}_{a} (F)+\D^{E}_{b} (B)$.
\end{enumerate}

\item If $\left|Y_N\right|  < b$, then there is a path 
$Q$ from $p(x)$ to $p(y)$ in $B$ which avoids $p(Y_N)$ of length $\ell(Q)\leq \D^{E}_{b-1} (B)\leq \D^{E}_{b} (B)$.
If $\left|E(F_x) \cap Y_D \right|\leq a$ or $\left|E(F_y) \cap Y \right|\leq a$, then obviously  $d_{G\setminus Y}(x,y) \leq \ell(Q) + \D^{E}_{a} (F)\leq \D^{E}_{b} (B)+ \D^{E}_{a} (F)$.\\
Now let  $\left|E(F_x) \cap Y \right|>a$ and $\left|E(F_y) \cap Y \right|>a$.
Then outside of the fibres $F_x$ and $F_y$ there are at most $b-1$ faulty edges.
Let $Y_N^{\prime}  \subseteq E(B)$ be set of edges from $p(x)$ to such neighbors $v_i\in V(B)$ for which fibre $F(v_i)$ contains faulty edges,
$Y_N^{\prime} =\{e=p(x)v\in E(B); \left|F(v) \cap Y_D \right|>0\}$.
Note that if $p(x)$ is a neighbor of $p(y)$ then  $p(x)p(y) \in Y_N^{\prime}$.
Then the subgraph $B \setminus (Y_N^{\prime} \cup p(Y_N))$ is  a connected graph as there are at most $b-1$ faulty edges outside of fibres $F_x$ and $F_y$.
Therefore there is a path $Q^{\prime}$ from $p(x)$ to $p(y)$ with $2\leq \ell(Q^{\prime})\leq \D^{E}_{b} (B)$ that avoids $p(Y_N)$ and there is no faulty edges in the fibre $F(v)$ of 
the  neighbor $v$ of $p(x)$ on the path $Q^{\prime}$, $e=p(x)v\subset Q^{\prime}$. 
Hence there is a path $x\rightarrow x^{\prime} \stackrel{P}{\rightarrow}y^{\prime}\stackrel{\tilde{Q}\setminus e}{\rightarrow} y$, which
avoids $Y$, where  $x^{\prime},y^{\prime} \in F(v)$ and $\ell(P)\leq \D(F)$. 
Thus $d_{G\setminus Y}(x,y)\leq 1+\D (F)+\ell(Q)-1 \leq \D^{E}_{a}(F)+\D^{E}_{b} (B)$.   \qed
\end{enumerate}
\end{itemize}
\end{pf*}

Clearly, Theorem \ref{EFD} improves Theorem \ref{EPsv} for all cases except in the following two  cases:
either  when $a=0$ and $F$ is a complete graph 
or  when  $b=0$ and $B$ is a complete graph.

\end{document}